\documentclass{amsart} 
\usepackage{amssymb}
\usepackage{amsmath}
\usepackage{amsfonts}

\begin{document}
\newtheorem{theo}{Theorem}[section]
\newtheorem{atheo}{Theorem*}
\newtheorem{prop}[theo]{Proposition}
\newtheorem{aprop}[atheo]{Proposition*}
\newtheorem{lemma}[theo]{Lemma}
\newtheorem{alemma}[atheo]{Lemma*}
\newtheorem{exam}[theo]{Example}
\newtheorem{coro}[theo]{Corollary}
\theoremstyle{definition}
\newtheorem{defi}[theo]{Definition}
\newtheorem{rem}[theo]{Remark}


\newcommand{\Bb}{{\bf B}}
\newcommand{\Cb}{{\mathbb C}}
\newcommand{\Nb}{{\mathbb N}}
\newcommand{\Qb}{{\mathbb Q}}
\newcommand{\Rb}{{\mathbb R}}
\newcommand{\Zb}{{\mathbb Z}}
\newcommand{\Ac}{{\mathcal A}}
\newcommand{\Bc}{{\mathcal B}}
\newcommand{\Cc}{{\mathcal C}}
\newcommand{\Dc}{{\mathcal D}}
\newcommand{\Fc}{{\mathcal F}}
\newcommand{\Ic}{{\mathcal I}}
\newcommand{\Jc}{{\mathcal J}}
\newcommand{\Kc}{{\mathcal K}}
\newcommand{\Lc}{{\mathcal L}}
\newcommand{\Oc}{{\mathcal O}}
\newcommand{\Pc}{{\mathcal P}}
\newcommand{\Sc}{{\mathcal S}}
\newcommand{\Tc}{{\mathcal T}}
\newcommand{\Uc}{{\mathcal U}}
\newcommand{\Vc}{{\mathcal V}}

\author{Nik Weaver}

\title [The liar paradox is a real problem]
       {The liar paradox is a real problem}

\address {Department of Mathematics\\
Washington University\\
Saint Louis, MO 63130}

\email {nweaver@math.wustl.edu}

\date{\em February 12, 2017}


\maketitle

\section{}

I recently wrote a book \cite{W} which, not to be falsely modest, I think
says some important things about the foundations of logic. So I
have been dismayed to see potential readers dismiss it after the first
few pages, which discuss the liar paradox. Their
probable opinion on that subject was eloquently expressed (not in reference
to my book) by the physicist Lubo\v{s} Motl \cite{M}:
\begin{quote}
Some fruitless discussion is dedicated to the Liar Paradox, i.e.\ even
more futile attempts to assign a truth value to the sentence
`this sentence is false'. This sentence can't be false because then
its negation would have to be right and the negation would imply that
the sentence is true; and vice versa. It can't be true and false at the
same moment which means that it can't be attributed any particular truth value.

I think that most intelligent schoolkids understand the previous paragraph.
On top of that, one may add some adults' interpretations. First of all,
there's nothing wrong about a proposition's having an undefined truth
value. Meaningless sentences can't be assigned truth values. That's the
case when the[y] fail to obey some rules of grammar or syntax. But even if
they do, they may fail to obey other conditions, conditions linked to the
`beef' of our axiomatic system.

This contradiction is avoided in any consistent framework to assign truth
values to some propositions because any such consistent framework does
forbid --- has to forbid (in order to be consistent) --- such sentences.
These requirements are reflected by special refinements such as the GB or
ZF set theory that overcome a similar paradox due to Bertrand Russell in
Georg Cantor's naive set theory (which allowed a set $M$ of all sets $X$
such that $X \not\in X$ which makes the question whether $M \in M$ equally
paradoxical) $\ldots$

The sequence of several paragraphs above really exhausts everything one
may say about the Liar Paradox $\ldots$
\end{quote}
I gather that this sort of view is quite widespread among mathematicians
and even logicians. So I felt it could be worthwhile to explain why it is
wrong. My intended audience is technically-minded people who think that
the liar paradox is trivially uninteresting. I understand that they come to
the subject already frustrated that anyone would even want to talk about it but
I believe I can give them clear reasons to reconsider their position.

(I should acknowledge right off that Kripke's `Jones versus Nixon' example
\cite{K} already convinced many philosophers of language to take these kinds
of paradoxes seriously. But I suppose it would not be as persuasive to a
technically-minded audience.)

\section{}

Motl's comment about ZFC --- the Zermelo-Frankel axiomatization
of set theory --- is a good place to start. He is right that it overcomes
the familiar paradoxes of naive set theory, by formalizing set-theoretic
reasoning in a way that accomplishes two goals:
\medskip

{\narrower{
\noindent 1. The paradoxes are blocked.
\medskip

\noindent 2. Ordinary unparadoxical reasoning is allowed.

\medskip}}

\noindent But Motl's vague reference to ``any consistent framework to assign
truth values to some propositions'' evades the following brute fact:
\medskip

{\narrower{
\noindent
{\underbar{We have no comparable axiomatization of truth.}}

\medskip}}

\noindent The vast technical literature on the subject contains many attempts,
but it has not yet produced a de facto standard formal system for reasoning
about truth that simultaneously blocks the liar paradox and retains ordinary
unparadoxical reasoning. Perhaps that in itself already gives some substance
to my claim that the liar paradox is a real problem.

Why has a workable formalization been so elusive?

\section{}

Is it due to incompetence? If so, then maybe there is an opportunity
here for someone unschooled in the subject, but with good sense, to put
his mind to the problem and clear everything up for us by formulating an
axiomatic system that does for truth what ZFC does for sets. It might be
easy!

I am being facetious, but I have come to suspect that many mathematicians
(of whom I am one) actually do believe incompetence to be the primary
explanation of our lack of a good formalization of truth. Perhaps they
know that philosophers think they can trap them in the following way:
\medskip

{\narrower{
\noindent P: Consider the sentence `This sentence is not true'. Is
it true or not?

\noindent M: It is neither true nor false, it has no truth value.

\noindent P: Aha, so in particular it is not true! But that is just
what it says of itself!

\medskip}}

\noindent and they have a rejoinder:
\medskip

{\narrower{
\noindent M: No, it does not `say' anything. It is meaningless.

\medskip}}

\noindent This simple comment apparently neutralizes any obstacle
the liar paradox may present to axiomatization. We can allow as valid a
formal version of the sentence `the liar sentence is neither true nor
false' (something Motl affirms), but block any subsequent inference of
the liar sentence itself by simply forbidding it as a possible assertion.
That should dissolve the paradox, or so one might think. As the logician
Arnon Avron puts it \cite{A},
\begin{quote}
Surely a meaningless sentence cannot say anything about
anything, in particular not about itself (or anything else). So
relying on what `it says of itself' depends on taking 
for granted that it is meaningful $\ldots$
I do wonder now if I am missing something here, and if so -
what can it possibly be.

Needless to say, for me the `liar sentences' of all types are indeed
completely meaningless, which is why I was never bothered by them $\ldots$
\end{quote}
What can he possibly be missing? Term substitution.

\section{}
I will explain. It is not difficult.

A sentence is just a string of symbols, so let us work with strings. We
can legitimately say things like
\begin{quote}
`sentence liar the' is not a true sentence.
\end{quote}
because `sentence liar the' is not even a grammatical sentence. Or,
equivalently, we can say
\begin{quote}
The concatenation of the strings `sentence li' and `ar the' is not
a true sentence.
\end{quote}
These two statements are trivially deducible from each other because the
{\it terms}
\begin{quote}
`sentence liar the'
\end{quote}
and
\begin{quote}
The concatenation of the strings `sentence li' and `ar the'
\end{quote}
evaluate to the same string. Agreed? The first term directly names the string
in question, while the second constructs it by a simple operation. Thus the
two sentences can be inferred from each other by a basic property of equality.
We can call this kind of inference {\it term substitution}.

An example of term substitution in a mathematical context: from
$$17 = 15 + 2$$
and
$$17\mbox{ is prime}$$
we can infer
$$15 + 2\mbox{ is prime}$$
--- in general, in any formal system that has an equality symbol, for
any constant terms $s$ and $t$ term substitution lets us derive the formula
$\phi(t)$ from the formula $\phi(s)$ plus the equality $s = t$.

Good. Now consider the following operation on strings. Given a string $S$,
replace every occurence of the symbol `{\#}' in $S$ with the full string
$S$ itself, in quotes. Call this operation {\it quining}. Thus, quining
`abc' just yields `abc', while quining `ab{\#}c' yields `ab`ab{\#}c'c'.
Then let $L$ be the following sentence:
\begin{quote}
The string obtained by quining `The string obtained by quining {\#} is not a
true sentence.' is not a true sentence.
\end{quote}
Is $L$ a true sentence?

Option 1: $L$ is a true sentence. A moment's thought shows that the string
it refers to is precisely itself. Since it is a true sentence, it does `say'
something, and what it says it that it is not true. Contradiction.

Option 2: $L$ is not a true sentence. The concern now is that this seems to
be just what $L$ says. But as Avron points out, that suggestion takes it for
granted that $L$ is meaningful. Can we agree that $L$ is not meaningful, that
it does not actually `say' anything?

Not really. Let $s$ be the term which simply quotes $L$ verbatim, and
let $t$ be the term
\begin{quote}
The string obtained by quining `The string obtained by quining {\#} is not a
true sentence.'
\end{quote}
Then $s$ and $t$ evaluate to the same string, namely the sentence $L$, so
we can go from
\begin{quote}
$s$ is not a true sentence.
\end{quote}
--- the premise of option 2 --- to
\begin{quote}
$t$ is not a true sentence.
\end{quote}
--- the paradoxical sentence $L$ itself --- by merely substituting $t$ for
$s$.
\medskip

{\narrower{
\noindent
{\underbar{Once you accept that $L$ is not a true sentence, you can immediately
infer}} {\underbar{$L$ by term substitution.}}

\medskip}}

Thus, in forbidding us from infering $L$ from the statement that $L$ is not
true, the Motl/Avron position commits itself to the conclusion that a
sentence of the form
\begin{quote}
$s$ is not a true sentence.
\end{quote}
is meaningful and true, but another sentence of the identical form
\begin{quote}
$t$ is not a true sentence.
\end{quote}
--- where $s$ and $t$ are concrete syntactic terms which evaluate to the
same string --- is meaningless.

To see how this would look in a formal setting, suppose ordinary Peano
arithmetic is augmented with a ``truth'' predicate symbol $T$ for
which we know $T[\bar{n}] \leftrightarrow \phi$ for any sentence $\phi$
with G\"odel number $n$. Then we
can produce a formal version $Q$ of the quining operation, which, given
the G\"odel number $n$ of a formula $\phi(x)$ with one free variable, returns
the G\"odel number of the formula $\phi(\bar{n})$. For a suitable choice
of $n$, the G\"odel number of the sentence
$$\neg T(Q[\bar{n}])$$
is $Q[\bar{n}]$ itself. This is just Tarski's
``truth'' version of G\"odel's construction of a sentence that denies its
own provability. Letting $k$ be the numerical value of $Q[\bar{n}]$, by
the truth property $T[\bar{k}]$ is equivalent to $\neg T(Q[\bar{n}])$,
and by term substitution the latter is equivalent to $\neg T[\bar{k}]$.
Which is a contradiction.

I suspect that anyone who has actually worked with formal theories
of truth would immediately reject the Motl/Avron idea as a non-starter,
for just this reason.

\section{}

I used the word `quining' for the construction discussed above because it is
similar to one introduced by W.\ V.\ O.\ Quine \cite{Q}. It is essentially a
natural language version of the construction employed by G\"odel in his
proof of the
first incompleteness theorem. (Incidentally, I learned about it from Douglas
Hofstadter's delightful book {\it G\"odel, Escher, Bach: An Eternal Golden
Braid} \cite{H}, which I read as an intelligent schoolkid. Did Motl, I wonder?)

I hope there is no question about the legitimacy of the quining operation.
There is nothing wrong with it: if strings are encoded numerically as
G\"odel numbers, it becomes a primitive recursive function.

Remember, I am trying to explain why it is so hard to formalize our
intuitive notion of truth in the same way that ZFC formalizes our intuitive
notion of sets. Any formal system for reasoning about strings with sufficient
number-theoretic resources to express primitive recursive functions is going
to have a quining operation. If it also has a truth predicate that can
be applied to its own sentences, and it allows term substitution, then
it cannot accomodate the Motl/Avron position.

\section{}

So as much as we want to say that liar sentences have no truth value, and
in particular are not true, if substitution of equivalent terms is a legal
deductive move, then this position leads to contradictions. The people who
take the liar paradox seriously are not just a bunch of fools. There are
real issues here which casual critics have not appreciated.

I suppose one could bite the bullet and affirm that mere term substitution
can sometimes turn a meaningful sentence into a meaningless one. But I would
be surprised to see anyone pursue this direction with much enthusiasm. You
then have to face up to the problem of how we can tell which instances of
term substitution are legitimate.

Kripke \cite{K} pointed out that in the natural language setting there
are many seemingly unremarkable sentences which are meaningful under most
circumstances, but which can become paradoxical if the empirical facts are
unfavorable. There are no empirical facts in pure mathematics, but essentially
the same problem is still present. If the language contains a truth predicate,
then for most numerical terms $t$ a sentence of
the form
\begin{quote}
The numerical value of $t$ is not the G\"odel number of a true sentence.
\end{quote}
may be straightforwardly true or false. But if the language
accomodates elementary number theory, then there will be cases
where, for example, $t$ evaluates to the G\"odel number of that very sentence,
rendering it paradoxical. It would not be hard for someone with knowledge of
elementary logic to design a primitive recursive term $t$ that contains
a numerical parameter $n$, and which evaluates to the G\"odel number of
the sentence `$2 + 2 = 4$' if $n$ is prime but which, if $n$ is composite,
evaluates to the G\"odel number of a sentence that affirms that the numerical
value of $t$ is not the G\"odel number of a true sentence. The moral is
that no simple syntactic filter is going to distinguish between ordinary and
paradoxical truth assertions.

\section{}

I do not mean to suggest that there are no other options; there are plenty.
One might hope to find a way out by ascribing truth not to sentences but
rather to the abstract ``propositions'' those sentences (allegedly) express.
Or by developing the idea that truth is situational, or indeterminate, or
subject to revision. Or more radically, by adopting a non-classical logic,
or even by accepting that some sentences can somehow simultaneously
be both true and not
true. In my experience, when a defender of the Motl/Avron view is shown the
term substitution problem, his first reaction is usually to make one of these
moves. Suffice it to say that they have all been thoroughly investigated, and
they all have serious problems.

One principal difficulty is that attempts to resolve the liar paradox often
turn out to involve concepts ``which, if allowed into the object language,
generate new paradoxes that cannot be dissolved by the account in question''
\cite{C}, a circumstance that is known in the liar literature as the ``revenge
problem''. For instance, in response to the idea that the truth of a sentence
may depend on the situation in which it is stated, we can formulate the
sentence ``This sentence is not true in any situation''. Sentence of this
sort are known in the literature as ``strengthened liars''.

In a word, my message to the intended audience for this paper is: whatever
simple idea you have for an easy resolution of the liar paradox --- we've
tried it, and it doesn't work.

\section{}

Surely the liar paradox is not important in itself. But to the extent that
it stands in the way of a satisfactory theory of truth, and that this is
a goal we care about, it matters.

Why should we care about having a satisfactory theory of truth? At this
point it may be helpful to look at the way truth ascriptions can appear
in mathematical reasoning, so we can see what is at stake.

By way of example, consider the formal system of Peano Arithmetic (PA).
Although the Peano axioms are quite simple, they are surprisingly strong
--- most ordinary number-theoretic reasoning can ultimately be reduced to
them. But not all number-theoretic reasoning: as G\"odel showed, if PA is
consistent then it cannot prove its own consistency. More precisely, it
cannot prove a certain number-theoretic sentence, ${\rm Con}({\rm PA})$,
which arithmetically expresses the consistency of PA.

Despite this, most of us think that PA is consistent. Why? How can we
be sure that there is no formal derivation of the sentence `$0=1$' in PA?
Because it has a model. If all of the Peano axioms are true in their
intended interpretation in the natural numbers, and if the inferential rules
of PA preserve truth, then we can argue by induction on the length of a
derivation that every theorem provable in PA is true in its intended
interpretation in $\mathbb{N}$. Since it is not true that zero equals
one in $\mathbb{N}$, the sentence `$0=1$' cannot be a theorem of PA.

This kind of reasoning requires a functioning concept of truth. If our
truth concept were inconsistent then such reasoning could not be trusted.

\section{}

However, this threat to Peano arithmetic might not be so serious, because it
is possible to explicitly define what we mean by truth in this particular
setting, and to thereby ensure that we have the tools we need to make the
consistency argument. That is, without referencing any general
concept of truth, we can define what it means for a formal sentence in
the language of arithmetic to be {\it true in $\mathbb{N}$}. The definition
goes by recursion on the length of the sentence. We start by saying under
what conditions on $s$ and $t$ an atomic sentence of the form $s = t$
is true in $\mathbb{N}$, where $s$ and $t$ are numerical terms, and then
we proceed inductively. ($A \wedge B$ is true in $\mathbb{N}$ if both
$A$ and $B$ are true in $\mathbb{N}$, etc.) This construction certainly
requires some infinitary resources in order to deal with quantifiers ---
indeed, more infinitary resources than are available within PA ---
but it does not require any prior concept of truth. And once we know what
it means for a formula of arithmetic to be true in $\mathbb{N}$, we
can make the argument for the consistency of PA sketched above.

Thus, we can quarantine the idea of arithmetical truth, give it a
formal definition, and put it to work without having to worry about
any liar related issues. This formal definition is rigorous and
philosophically unproblematic (at least, modulo concerns about
infinity, which are a separate matter). Moreover, it generalizes to
any formal language with a set model, a result due to Tarski and Vaught
\cite{TV}.

\section{}

We have made definite progress. The pessimism of our earlier discussion,
with its implication that an axiomatic theory of truth cannot hope to
satisfy even the humblest conditions we might impose, has been countered
by a definite theorem which effectively states that in a wide variety of
settings a rigorous definition of truth is available. 

The contrast is striking: when we try to analyze truth as a universal concept,
we just get one contradiction after another, but when we localize to any
specific formal language that has a well-defined set-theoretic interpretation,
we can give an explicit, almost trivial construction of something that we
intuitively recognize as a truth predicate. The natural conclusion is that
the global concept is fictional; all we have are local truth predicates.

Another way to say this is that we should be thinking not in terms of a
concept of `truth', but in terms of a concept of `truth predicate'. This was
Tarski's view. That is, the right question to ask is not `What is truth?',
it is `What counts as a truth predicate for a given language?' And he has
an answer to this question \cite{T}. According to Tarski, the functional
role of truth is encapsulated in the {\it T-scheme}
$$T(\mbox{`$\mathcal{A}$'}) \leftrightarrow \mathcal{A},$$
where $T(\cdot)$ is a candidate truth predicate and $\mathcal{A}$ is to be
replaced by any syntactically correct sentence. The classic example is
`$\,$`Snow is white' is true if and only if snow is white'. Tarski's idea
is to use the T-scheme, with $\mathcal{A}$ ranging over all the sentences
of some language, to say what it means for a predicate $T$ (in formal
terms, a formula with one free variable) to count as a truth predicate
for that language. This predicate would itself typically belong not to
the target language, but to some metalanguage.

It may look as though the T-scheme trivializes truth, if all it tells
us is that asserting a sentence is the same as asserting its truth. But
this equivalence gives us the ability to, for example, say something
about snow (that it is white) by saying something about the syntactic
string `snow is white' (that it is true), and depending on our ability to
reason about and manipulate syntactic strings, this can be a powerful tool.
Most importantly, we can quantify over syntactic strings. Recall
the argument for the consistency of Peano arithmetic: it relied on general
principles like `For any sentences $A$ and $B$, if $A$ and $A \to B$ are
both true, then so is $B$.' One needs a truth predicate to make statements
like this, and one needs the general fact, not any particular instance, to
make the inductive argument for consistency.

\section{}

Tarski's account seems to remove the mystery surrounding truth. On this view
truth is not some occult metaphysical quality; when we talk about `truth' we
are merely talking about any predicate that makes the T-scheme work. The
Tarski-Vaught theorem shows that we can construct truth predicates in a
variety of formal settings, and the lesson of the liar paradox is that this is
the best we can hope for: nothing can function as a truth predicate globally.
End of story.

Or is it? Many commentators have been uncomfortable with the idea that there
is no single, universal concept of truth. It seems wrong. The metalanguage
in which a Tarskian truth predicate appears can itself be the target of a
new, broader truth predicate, leading to the idea of hierarchies of truth
predicates of ever-increasing generality. Again, this violates our raw
intuition of truth as a unitary quality.

A more focussed objection is that the T-scheme is, in fact, disastrously
circular. Just above I said, with intentional vagueness, that a truth
predicate is one which makes the T-scheme `work'. What this means is,
precisely, that it must make every instance of the T-scheme
{\underbar{true}}. Yes. Affirming that the T-scheme holds for a given
predicate $T$ requires the use of a preexisting truth predicate, just
as much as affirming the law `if $A$ and $A \to B$ are both true, then
so is $B$' does. Any single instance
of this law can be expressed without using truth, just by substituting
specific sentences in the template `if $A$ and $A \to B$ then $B$'. But
quantifying over all $A$ and $B$ requires a truth predicate. The T-scheme
is in exactly the same situation. Any single instance can be expressed
without using a truth predicate, but you need a truth predicate to express
that something happens for every sentence $A$. That is just the kind of
thing truth predicates are good for.

This is not a pedantic complaint. In order to use the T-scheme to say what
counts as a truth predicate, you need to use a truth predicate. Which is
to say:
\medskip

{\narrower{
\noindent
{\underbar{We cannot use the $T$-scheme to say what truth means without
already}} {\underbar{knowing what truth means.}}

\medskip}}

\noindent More precisely, we cannot use it to say what truth means in a
given context without already knowing what it means in a broader context.
I call this the {\it Tarskian catastrophe}. It shows that we have not, after
all, succeeded in eliminating truth as a primitive notion. The problem is
easy to miss because we are so used to using truth to convert schemes into
sentences that we may not even notice the need to use it in that way here.

\section{}

Tarski was aware of this difficulty, and his response was to reframe the
condition in terms of each instance of the scheme not being true, but being
formally provable in some metasystem. In effect, we take all the separate
instances of the T-scheme as a collective definition, an axiom scheme. But
this completely vitiates the idea of truth. There is an essential difference
between proving each instance of a scheme and proving a single statement
which implies every instance. If all we have is the T-scheme as a scheme,
then we have no basis to affirm that for all $A$ and $B$, if $A$
and $A \to B$ are both true, then so is $B$. The negation of this general
law is consistent with any finitely many instances of the T-scheme, and
therefore the general law is not a logical consequence of the T-scheme
as a scheme.

The problem is fundamental because practically the whole point of having a
truth predicate in mathematical settings is to enable ourselves to make
general statements about truth. If the T-scheme is verified in a
piecemeal way then it can only be used in a piecemeal way, which is to
say, it cannot be used substantively.

In order to characterize truth using the T-scheme we have to find a way to
globally affirm every instance of the scheme, something Tarski's dodge fails
to accomplish. This means that the liar paradox is still a problem, because
it turns out that we still require
a global notion of truth, and that still runs afoul of the paradox.

\section{}

One might ask why we even need the T-scheme when we already have the
Tarski-Vaught theorem which explicitly constructs a truth predicate
for any interpreted language. Why not just use this as a definition?

Because the Tarski-Vaught theorem does not actually do this. As stated, it
only applies to languages equipped with a set model, i.e., languages for
which we have a set which the variables of the language are understood
as ranging over. (For instance, in the case of PA the variables are
supposed to range over the set $\mathbb{N}$.) But in mathematics we often
want to make assertions about variables which range over a proper class,
a situation that would not be covered. One may be willing to dismiss class
models as unnecessary extravagance, but if we want to use the Tarski-Vaught
theorem as a definition of `truth predicate' that option is not available.
This theorem itself refers to arbitrary set models, and hence requires for
its own expression a language with variables that range over all sets.
This is a necessity, not a luxury.

It would not be hard to generalize the theorem to apply to
languages equipped with class models, but in order to do this one would
need to use a language with variables that range over all classes, a language
to which the theorem still would not apply. Indeed, however broadly we
generalize the Tarski-Vaught construction, we can never get it to apply
to the very language in which the generalization itself is effected. For
then this language would be able to supply its own truth definition $\ldots$
leading to a liar paradox.

\section{}

I know how to solve these problems \cite{W}. My idea is to segregate from
truth a distinct but related notion, the notion of assertibility. Informally,
a sentence is {\it assertible} if we have an unimpeachable right to affirm
it. This is something intuitionists have been very interested in, but
their analysis of it seems wrong to me. (They want it to replace the
classical notion of truth, which they reject, but they also think it
satisfies the T-scheme, which would effectively make it a classical
truth predicate.) I do think they are right on one crucial point, that
reasoning about assertibility demands the use of intuitionistic logic.

This is not the place for a detailed explanation, but to summarize some
of my conclusions, I find that assertibility can be given a simple,
elegant axiomatization, one that is (at least in various limited settings)
provably consistent. The assertibility version of the liar paradox is
defeated by rejecting some instances of the assertibility T-scheme ---
a subtle point, because it hinges on an intuitionistic interpretation of
implication. It is assertibility, or constructive truth,
which is, in my view, the global concept, while classical truth is,
just as Tarski said, only something which makes sense locally. The
explanation for our naive intuition of truth as a unitary quality is that
it is actually an intuition for assertibility, or maybe that it involves
some equivocation between classical and constructive truth. The Tarskian
catastrophe is resolved by saying that $T$ is a truth predicate for an
interpreted language if all instances of the T-scheme, with $A$
ranging over the sentences of the language, are assertible. Indeed,
since it is truly global, the concept of assertibility is available to
handle any number of other tasks which classical truth, because of its
inherently local nature, cannot. A successful treatment of second order
logic, for instance --- something which has been a problem since Frege's
approach was thwarted by Russell's paradox --- now becomes possible.
In short, I find that the constructive notion of truth can be used to
do all the things one wants classical Tarskian truth to do that it cannot.

This is a grand claim, but I am confident of it. However, before one
can even start to understand such things, one has to get past this inane
idea that the liar paradox is not a real problem. It absolutely is.


\end{document}